\documentclass[12pt]{article}
\usepackage[reqno]{amsmath}
\usepackage{amssymb, amsthm, enumerate}
\usepackage{graphicx}
\usepackage{tikz}
\usetikzlibrary{decorations.markings}
\usepackage{array,booktabs,hhline}
\usepackage{float}
\usepackage{wrapfig}
\usepackage{url,hyperref}
\usepackage{longtable}

\expandafter\ifx\csname pdfpagebox\endcsname\relax\else
\fi

\makeatletter
\newtheorem*{rep@theorem}{\rep@title}
\newcommand{\newreptheorem}[2]{%
\newenvironment{rep#1}[1]{%
 \def\rep@title{#2 \ref{##1}}%
 \begin{rep@theorem}}%
 {\end{rep@theorem}}}
\makeatother

\newtheoremstyle{plainsl}%
	{\topsep}
	{\topsep}
	{\slshape} 
	{}
	{\normalfont\bfseries}
	{.}
	{ }
	{}

\swapnumbers

{\theoremstyle{plainsl}
\newtheorem{theorem}{Theorem}[section]
\newtheorem{lemma}[theorem]{Lemma}
\newtheorem{corollary}[theorem]{Corollary}}
{\theoremstyle{remark}
}

\newreptheorem{theorem}{Theorem}

\newcommand\cref[1]{Corollary~\ref{cor:#1}}

\def\proof{\noindent{{\sl Proof. }}}
\def\sqr#1#2{{\vbox{\hrule height.#2pt
    \hbox{\vrule width.#2pt height#1pt \kern#1pt
        \vrule width.#2pt}\hrule height.#2pt}}}
\def\eqed{\sqr53}
\def\qed{%
    \ifmmode\eqno\eqed
    \else\nobreak\ \hfill\eqed\medbreak\fi}

\newcommand\cA{{\mathcal A}}

\newcommand\cD{{\mathcal D}}

\newcommand\cx{{\mathbb C}}

\newcommand\re{{\mathbb R}}

\newcommand\comp[1]{{\mkern2mu\overline{\mkern-2mu#1}}}

\newcommand\seq[3]{#1_{#2},\ldots,#1_{#3}}

\newcommand\cmmnt[1]{}
\newcommand\ip[2]{\langle#1,#2\rangle}
\DeclareMathOperator{\rk}{rk}
\DeclareMathOperator{\tr}{tr}

\newcommand\mat[3]{\mathrm{Mat}_{#1\times #2}(#3)}
\newcommand\matnn[1]{\mat{n}{n}{#1}}
\newcommand\ones{{\mathbf 1}}
\DeclareMathOperator{\im}{im}

\DeclareMathOperator{\cmm}{Comm}
\DeclareMathOperator{\vc}{vec}

\renewcommand\th{\theta}
\newcommand\mxm{\widehat{M}}

\title{A New Perspective on the \\ Average Mixing Matrix}

\author{Gabriel\! Coutinho\thanks{Universidade Federal de Minas Gerais, Belo Horizonte, Brazil. \url{gabriel@dcc.ufmg.br}}, Chris\! Godsil\thanks{University of Waterloo, Waterloo, Canada. \url{{cgodsil,
h3zhan}@uwaterloo.ca}}, Krystal\! Guo\thanks{Universit\'{e} libre de Bruxelles, Brussels, Belgium. \url{guo.krystal@gmail.com}} \thanks{This work was done when K. Guo was a postdoctoral fellow at the University of Waterloo.}, Hanmeng\! Zhan\footnotemark[2]}

\begin{document}
\maketitle

\begin{abstract}
	We consider the continuous-time quantum walk defined on the adjacency matrix of a graph. At each instant, the walk defines a mixing matrix which is doubly-stochastic. The average of the mixing matrices contains relevant information about the quantum walk and about the graph. We show that it is the matrix of transformation of the orthogonal projection onto the commutant algebra of the adjacency matrix, restricted to diagonal matrices. Using this formulation of the average mixing matrix,  we find connections between its rank and  automorphisms of the graph.
\end{abstract}

\section{Introduction}\label{sec:intro}

Let $X$ be a graph. We are interested in \textsl{continuous quantum walks} on $X$, which we will define
now. The states of the walk are represented by density matrices, positive semidefinite matrices
with rows and columns indexed by the vertex set $V(X)$ of $X$, and having trace $1$. If $A$ is the 
adjacency matrix of $X$, we define the \textsl{transition matrix} $U(t)$ of the walk by
\[
	U(t) = \exp(itA).
\]
 The continuous-time quantum walk is an important object of study in quantum computing since it is an universal computational primitive \cite{ChildsUniversalQComputation}. They were first studied in \cite{FarhiGutmann}. Since then, many aspects of quantum walks have been studied, including state transfer \cite{KayReviewPST,GodsilStateTransfer12,CoutinhoGodsilGuoVanhove2,VinetZhedanovHowTo}, uniform mixing \cite{AdaChanComplexHadamardIUMPST,GodsilMullinRoy,TamonAdamczakUniformMixingCycles} and average mixing \cite{GodsilAverageMixing,TamonAdamzackAverageMixing}. 
We note that $U(t)$ is unitary and symmetric. If the initial state of the system is given by a 
density matrix $D$, then the state of the system at time $t$ is
\[
	U(t)DU(-t).
\]
We note that we do not have direct access to states, all that experiment provides is the value
of expressions of the form
\[
	\tr(DP_r)
\]
where $\seq P1n$ are positive semidefinite matrices such that $\sum_rP_r=I$. Because of this
we find ourselves using the \textsl{trace inner product} on $\mat mn\cx$:
\[
	\ip{M}{N} = \tr(M^TN).
\]

Unlike classical random walks on a connected graph, a continuous quantum walk does not reach a steady
state. Thus, if the eigenvalues of $X$ are integers (e.g., if $X$ is the complete graph $K_n$)
then $U(t)$ is a periodic function of $t$. The focus of this paper, the \textsl{mixing matrix} $\mxm$
of $X$ provides a useful substitute for a steady state.

There are two convenient ways to define $\mxm$. The first is to introduce what we call the
mixing matrices $M(t)$ of the walk, given by
\[
	M(t) = U(t) \circ \comp{U(t)}.
\]
(Here we use $M\circ N$ to denote the Schur or element-wise product of two matrices of the same order.)
Since $U(t)$ is unitary and symmetric, we see that $M(t)$ is symmetric non-negative matrix
with each row and column summing to one. We can now define
\begin{equation}
	\label{eq:mxmTto}
	\mxm = \lim_{T\to\infty}\frac1T \int_0^T M(t)\,dt.
\end{equation}
To see that this definition makes sense, we recall that since $A$ is real and symmetric, it
has a spectral decomposition
\[
	A = \sum_r \th_r E_r
\]
where $\th_r$ runs over the distinct eigenvalues of $A$ and $E_r$ is the matrix that represents
orthogonal projection onto the $\th_r$-eigenspace of $A$. Given this we also have
\[
	U(t) = \sum_r e^{it\th_r}E_r
\]
and, as $\comp{U(t)}=U(-t)$,
\[
	M(t) = \sum_{r,s} e^{it(\th_r-\th_s)} E_r\circ E_s.
\]
It follows that the limit on the right side of \eqref{eq:mxmTto} is exists and is equal to
\[
	\sum_r E_r^{\circ2},
\]
that is, the limit is the sum of the Schur squares of the spectral idempotents of $A$.

We consider one example. The spectral idempotents of the complete graph $K_n$ are
\[
	\frac1n J,\quad I-\frac1n J,
\]
(with corresponding eigenvalues $n-1$ and $-1$) and hence its average mixing matrix is
\[
	\left(1-\frac2n\right)I + \frac2{n^2}J.
\]
For large $n$, this is very close to $I$. This is unexpected, but we are dealing with quantum 
physics, where the unexpected is not uncommon.

If $A\in\mat nn\re$, its  \textsl{commutant} $\cmm(A)$ is the set of matrices that commute with $A$.
This is an object of some combinatorial interest---the permutation matrices in $\cmm(A)$ are
the automorphisms of $A$. The first main result of this paper is that, relative to a natural
basis, $\mxm$ is the matrix that represents the restriction to the diagonal matrices of
the orthogonal projection from $\mat nn\re$ onto $\cmm(A)$. We then use this connection
to investigate the relations between properties of the graph $X$ and the rank of $\mxm$.

\section{Projection onto the commutant of the adjacency matrix}

Let $A$ be a $n\times n$ real symmetric matrix with columns and rows indexed by elements of a set $V$. We denote by $\mat nn\re$ the set of $n\times n$ matrices with entries in $\re$. Let the spectral decomposition of $A$ be given as follows:
\[
	\sum_{r=1}^d \theta_r E_r,
\]
where $\seq\th1m$ are the distinct eigenvalues of $A$ and $E_r$ is the projection matrix onto the $\theta_r$ eigenspace of $A$.  We denote by $\cmm(A)$ the set of all real matrices which commute with $A$. 

\begin{lemma}\label{lem:dimcomm}
	If the eigenvalues $\seq\th1d$ of $A$ have multiplicities $\seq m1d$ respectively,
	them the dimension of $\cmm(A)$ is $\sum_r m_r^2$.
\end{lemma}

\proof
When $A$ is diagonal, this is immediate. Since $A$ is similar to a diagonal matrix, 
the lemma follows.\qed

We consider a map  $\Psi: \re_{n\times n} \rightarrow \cmm(A)$ such that 
\[
	\Psi(M) = \sum_{r=0}^d E_r M E_r,
\]
for $M \in \re_{n\times n}$. We see that $\Psi$ is an endomorphism of $\re_{n\times n}$. 

\begin{lemma}\label{lem:psilem} 
	The following are true:
	\begin{enumerate}[(i)] 
		\item $\Psi$ is idempotent. 
		\item The image of $\Psi$ is equal to $\cmm(A)$. 
		\item $\Psi$ is self-adjoint; that is $\ip{M}{\Psi(N)} = \ip{\Psi(M)}{N}$.
	\end{enumerate}
\end{lemma} 

\proof 
Since $E_r^2=E_r$ and $E_rE_s=0$ if $r\ne s$, it is immediate that $\Psi^2(M)=\Psi(M)$
for any matrix $M$. As 
\[
	A E_r M E_r = \th_r E_r M E_r = E_r M E_r A
\]
we see that $\Psi(M)\in\cmm(A)$ for any $M$. Each idempotent $E_r$ is a polynomial in $A$
and consequently if $N\cmm(A)$, then
\[
	\Psi(N) = \sum_r E_rNE_r = \sum_r NE_r^2 = \sum_r NE_r = N\sum_r E_r = NI = N.
\]
Hence each element of $\cmm(A)$ lies in the image of $\Psi$ and therefore $\im(\Psi)=\cmm(A)$.

Finally
\[
	\ip{M}{E_rNE_r} = \tr(M^TE_rNE_r) = \tr(E_rM^TE_rN) = \ip{E_rME_r}{N},
\]
from which it follows that $\Psi$ is self-adjoint.\qed

\begin{corollary} 
	The map $\Psi$ is the orthogonal projection of the $n\times n$ real matrices onto $\cmm(A)$. 
\end{corollary} 

\proof 
This is immediate from the fact that $\Psi$ is idempotent and self-adjoint---but we simply note
that for any two $n\times n$ matrices $M$ and $N$
\[
	\ip{\Psi(M)}{N-\Psi(N)} = \ip{M}{\Psi(N-\Psi(N))} =\ip{M}{\Psi(N)-\Psi^2(N)} = 0.\qed
\]

The following is a standard fact in linear algebra, see e.g., \cite[Lemma 4.3.1]{HornJohnsonTopics}. 

\begin{lemma} 
	For $B,C,N \in \re_{n\times n}$, we have that 
	\[
		\vc(C N B^T) = (B\otimes C) \vc(N).\qed
	\] 
\end{lemma}

From this we see that, relative to the standard basis of $\matnn\re$, the matrix that represents 
$\Psi$ is $\sum_r E_r \otimes E_r$. (Since $\tr(E_r)=m_r$ it follows that $\tr(\Psi)=\sum_r m_r^2$;
since $P$ is idempotent $\tr(P)=\rk(P)$ and thus we have a second proof of Lemma~\ref{lem:dimcomm}.)

\section{Average states}

Recall that a state $D$ is a positive semidefinite matrix with trace $1$. If $D$ is the initial
state of a continuous quantum walk, then the state $D(t)$ at time $t$ is
\[
	D(t) = U(t)DU(-t).
\]
Using the spectral decomposition of $U(t)$, we have
\[
	D(t) = \sum_{r,s} e^{it(\th_r-\th_s)}E_r D E_s
\]
whence
\[
	\lim_{T\to\infty}\frac1T \int_0^T D(t)\,dt = \sum_r E_rDE_r.
\]
We call this limit (or sum) an \textsl{average state} and denote it by $\Psi(D)$. We see at once 
that $\Psi(D)$ is equal to the orthogonal projection of $D$ onto $\cmm(A)$.

Since $D$ is positive semidefinite and the idempotents $E_r$ are symmetric, we see that $\Psi(D)$
is a positive semidefinite matrix. As
\[
	\tr(\Psi(D)) = \sum_r\tr(E_rDE_r) = \sum_r \tr(DE_r^2) = \sum_r \tr(DE_r) = \tr(D)
\]
we also see that $\Psi(D)$ is a density matrix.

In the context of quantum walks on graphs, there is a natural class of density matrices we
will focus on. If $a\in V(X)$, let $e_a$ denote the standard basis vector of $\cx^{V(X)}$
indexed by $a$ and define $D_a=e_ae_T$. Then certainly $D_a$ is a density matrix, moreover
$\rk(D_a)=1$ and $D_a^2=D_a$. (Physicists refer to a density matrix with rank $1$ as
a \textsl{pure state}.)

The following theorem gives one reason why average states are of interest.

\begin{theorem}\label{thm:avmxGram}
	Let $X$ be a graph. The average mixing matrix is the Gram matrix of the average states
	$\Psi(D_a)$ for $a$ in $V(X)$.
\end{theorem}

\proof
Our claim is that, if $a,b\in V(X)$, then
\[
	\ip{\Psi(D_a)}{\Psi(D_b)} = (\mxm)_{a,b}.
\]
Now
\[
	\Psi(D_a)\Psi(D_b) = \sum_r E_rD_aE_r \> \sum_s E_sD_bE_s
\]
and so
\[
	\tr(\Psi(D_a)\Psi(D_b)) = \sum_r \tr(E_rD_aE_rD_bE_r).
\]
Further
\[
	E_rD_aE_rD_bE_r = E_r e_ae_a^TE_re_be_b^T E_r = (E_r)_{a,b}\, E_r e_ae_b^T E_r
\]
and therefore
\[
	\tr(\Psi(D_a)\Psi(D_b)) = \sum_r(E_r)_{a,b}(E_r)_{b,a} = \sum_r(E_r^{\circ2})_{a,b}
\]
and the theorem follows.\qed

\begin{corollary}
	The dimension of the space spanned by the average state $\Psi(D_a)$ for $a\in V(X)$ is
	equal to $\rk(\mxm)$.\qed
\end{corollary}

\begin{corollary}
	The average mixing matrix is completely positive semidefinite.\qed
\end{corollary}

\section{Diagonal matrices in the commutant}

In this section, we consider the restriction of $\Psi$ to the set of $n\times n$ diagonal matrices. 
Let $\cD$ be the set of $n\times n$ diagonal matrices. The standard basis of $\cD$ is given 
by $\{D_a\}_{a\in V}$ where $D_a := e_ae_a^T$.

The map that sends a matrix $M$ to $\Psi(M)\circ I$ is a linear map from $\matnn\re$
into the space $\cD$ of diagonal matrices. We denote the restriction of this map to $\cD$
by $\Phi$; it is evidently an endomorphism of $\cD$.

\begin{lemma} \label{lem:rkM}
	Relative to the standard basis $\{D_a\}_{a\in V}$ of $\cD$, the matrix that represents
	$\Phi$ is $\mxm$. Hence $\rk(\mxm)=\dim(\Phi(\cD))$.
\end{lemma}

\proof
The entries of the matrix representing $\Psi$ are given by the the inner 
products $\ip{D_a}{\Psi(D_b)}$ for vertices $a$ and $b$ of $X$. We have
\[
	\ip{D_a}{\Psi(D_b)} = \sum_r \tr(D_aE_rD_bE_r)
\]
and
\[
	\tr(D_aE_rD_bE_r) = \tr(e_a e_a^TE_re_b e_b^TE_r) = (E_r)_{a,b}(E_r)_{b,a} 
		= (E_r^{\circ2})_{a,b}.
\]
Hence
\[
	\ip{D_a}{\Psi(D_b)} = \mxm_{a,b}.\qed
\]

\begin{lemma}
	If $D$ is diagonal, then $\Psi(D)=0$ if and only if $I\circ\Psi(D)=0$.
\end{lemma}

\proof
Since $\ker(\Psi)=\cmm(A)^\perp$, we see that $\Psi(D)=0$ if and only if $M\in\cmm(D)^\perp$.
Note that if $D$ is diagonal, then
\[
	\ip{D}{N} = \ip{D}{I\circ N}
\]

Assume $D$ is diagonal and $I\circ\Psi(D)=0$. Then for each vertex $a$ of $X$, we have
\[
	0 = \ip{D_a}{I\circ\Psi(D)} = \ip{D_a}{\Psi(D)} = \ip{\Psi(D_a)}{\Psi(D)}
\]
and so $\Psi(D)$ is orthogonal to each matrix $\Psi(D_a)$. Since $\Psi(D)$ lies in the span
of the matrices $\Psi(D_a)$, we conclude that $\Psi(D)=0$.\qed

Let $\cA_0$ denote the set of matrices in $\cmm(A)$ with all diagonal entries equal to $0$; that is 
\[
	\cA_0 = \{N \in \cmm(A): N\circ I = 0\}.
\]
Observe that $\cA_0$ is always non-empty, since $A \in A_0$. Also $I\in\Psi(D)$, and therefore
the direct sum decomposition in our next result is always non-trivial.

\begin{lemma} \label{lem:decomp}
	$\cmm(A) = \Psi(\cD) \oplus \cA_0$. 
\end{lemma}

\proof
If $M\in\cA_0$, then $\Psi(M)=M$ and so
\[
	\ip{\Psi(D_a)}{M}  = \ip{D_a}{\Psi(M)} = \ip{D_a}{M} = \tr(e_ae_a^TM) = M_{a,a}.
\]
Accordingly $M$ is orthogonal to each matrix $\Psi(D_a)$ for $a$ in $V(X)$ if and only
if $M\circ I=0$.\qed

\section{The rank of some average mixing matrices}

There are a number of graph invariants that can be constructed from the 
average mixing matrix. In this section we focus on $\rk(\mxm)$. From
\cite{GodsilAverageMixing} we know that if $\rk(\mxm)=1$, then $X$ is $K_1$ or $K_2$.

Our first two results concern graphs with only simple eigenvalues; we note that 
almost all graphs have this property \cite{TaoVanVuRandomSimpleSpectrum}.

\begin{lemma}\label{lem:simps}
	Assume $X$ is a graph with simple eigenvalues on $n$ vertices. If $n\ge2$, 
	then $\rk(\mxm) \leq n-1$. Further, if $X$ is regular and $n\ge4$, then 
	$\rk(\mxm)\leq n-3$.
\end{lemma}

\proof
Since the eigenvalues of $X$ are simple, $\dim(\cmm(A))=n$ and since $A\in\cA_0$,
it follows from Lemma~\ref{lem:decomp} that $\dim(\Psi(\cD))\le n-1$.

Suppose $X$ is regular with valency $k$. Then the matrices $A$, $J-I$ and $A^2-kI$
all lie in $\cA_0$. If these matrices are linearly dependent, the minimal polynomial
of $A$ has degree at most two and hence $A$ has at most two eigenvalues. 
Since $X$ has $n$ distinct eigenvalues, the minimal polynomial of $A$ over 
$\ones^{\perp}$ has degree $n-1$. Thus $n-1 \leq 2$, which contradicts our 
assumption that $n\ge4$.\qed

Table~\ref{tab:simpevs} shows the number of graphs on $n$ vertices whose mixing matrices have given ranks. The only regular graph with simple eigenvalues on up to $8$ vertices is $K_2$, so the data for regular graphs is not included in this table. Instead, we count cubic graphs with simple eigenvalues on $10$ to $18$ vertices, as presented in Table~\ref{tab:cubicsimpevs}.

\begin{table}[htbp]
    \centering
    \begin{tabular}{|c|c|c|c|}
    \hline \rule{0pt}{2.8ex} 
    $n$ & $\rk(\mxm)$ & \# graphs & \# simple eigenvalues\\
    \hhline{|=|=|=|=|}
    3&  2&  1&  1\\
    \hline
    3&  3&  1&  0\\
    \hhline{|=|=|=|=|}
    4&  2&  3&  2\\
    \hline
    4&  3&	1&  1\\
    \hline
    4&	4&	2&	0\\
    \hhline{|=|=|=|=|}
    5&	3&	11&	8\\
    \hline
    5&	4&	6&	3\\
    \hline
    5&	5&	4&	0\\
    \hhline{|=|=|=|=|}
    6&	2&	2&	2\\
    \hline
    6&	3&	27&	12\\
    \hline
    6&	4&	32&	21\\
    \hline
    6&	5&	35&	19\\
    \hline
    6&	6&	16&	0\\
    \hhline{|=|=|=|=|}
    7&	3&	6&	5\\
    \hline
    7&	4&	189&	121\\
    \hline
    7&	5&	240&	158\\
    \hline
    7&	6&	352&	255\\
    \hline
    7&	7&	66&	0\\
    \hhline{|=|=|=|=|}
    8&	2&	3&	3\\
    \hline
    8&	3&	39&	25\\
    \hline
    8&	4&	466&	236\\
    \hline
    8&	5&	1360&	776\\
    \hline
    8&	6&	2523&	1492\\
    \hline
    8&	7&	5781&	4787\\
    \hline
    8&	8&	945&	0\\
    \hline
    \end{tabular}
    \caption{Number of graphs with given $\rk(\mxm)$.}
    \label{tab:simpevs}
\end{table}

  \begin{table}[p]
  \begin{tabular}{|c|c|c|c|}
    \hline \rule{0pt}{2.8ex} 
    $n$ & $\rk(\mxm)$ & \# cubic  & \# simple \\
    & &graphs &eigenvalues\\
    \hhline{|=|=|=|=|}
    10&	3&	2&	2\\
    \hline
    10&	5&	8&	1\\
    \hline
    10&	6&	5&	3\\
    \hline
    10&	7&	1&	0\\
    \hline
    10&	10&	3&	0\\
    \hhline{|=|=|=|=|}
    12&	3&	1&	0\\
    \hline
    12&	4&	3&	0\\
    \hline
    12&	5&	8&	3\\
    \hline
    12&	6&	11&	2\\
    \hline
    12&	7&	18&	6\\
    \hline
    12&	8&	14&	4\\
    \hline
    12&	9&	14&	3\\
    \hline
    12&	10&	11&	0\\
    \hline
    12&	11&	2&	0\\
    \hline
    12&	12&	3&	0\\
    \hhline{|=|=|=|=|}
    14&	4&	13&	12\\
    \hline
    14&	5&	19&	12\\
    \hline
    14&	6&	30&	7\\
    \hline
    14&	7&	82&	37\\
    \hline
    14&	8&	97&	65\\
    \hline
    14&	9&	66&	37\\
    \hline
    14&	10&	62&	45\\
    \hline
    14&	11&	117&	101\\
    \hline
    14&	12&	18&	0\\
    \hline
    14&	13&	3&	0\\
    \hline
    14&	14&	2&	0\\
    \hline
    \end{tabular} 
      \begin{tabular}{|c|c|c|c|}
      \hline \rule{0pt}{2.8ex} 
    $n$ & $\rk(\mxm)$ & \# cubic  & \# simple \\
    & &graphs &eigenvalues\\
    \hhline{|=|=|=|=|}
    16&	4&	4&	3\\
    \hline
    16&	5&	45&	29\\
    \hline
    16&	6&	58&	11\\
    \hline
    16&	7&	122&	49\\
    \hline
    16&	8&	252&	112\\
    \hline
    16&	9&	393&	220\\
    \hline
    16&	10&	359&	144\\
    \hline
    16&	11&	311&	141\\
    \hline
    16&	12&	684&	464\\
    \hline
    16&	13&	1365&	1008\\
    \hline
    16&	14&	366&	0\\
    \hline
    16&	15&	77&	0\\
    \hline
    16&	16&	24&	0\\
    \hhline{|=|=|=|=|}
    18&	5&	48&	45\\
    \hline
    18&	6&	147&	59\\
    \hline
    18&	7&	226&	78\\
    \hline
    18&	8&	414&	152\\
    \hline
    18&	9&	1268&	724\\
    \hline
    18&	10&	1785&	982\\
    \hline
    18&	11&	1865&	842\\
    \hline
    18&	12&	1264&	539\\
    \hline
    18&	13&	1940&	1146\\
    \hline
    18&	14&	7254&	5819\\
    \hline
    18&	15&	19302&	16060\\
    \hline
    18&	16&	4763&	0\\
    \hline
    18&	17&	643&	0\\
    \hline
    18&	18&	382&	0\\
    \hline
     \end{tabular} 
    \caption{Number of cubic graphs with given $\rk(\mxm)$.}
    \label{tab:cubicsimpevs}
    \end{table}

\begin{lemma}
	Assume $X$ is a graph on $n$ vertices with only simple eigenvalues.
	If $X$ is bipartite, then $\rk(\mxm) \le \lfloor(n+1)/2\rfloor$.
\end{lemma}

\proof
Note that $A^k \in \cA_0$ if and only if $k$ is odd, and for all $k < n$, 
these matrices are independent. Hence $\dim \cA_0 \ge \lfloor n/2\rfloor$, and the
result follows from Lemmas~\ref{lem:rkM} and \ref{lem:decomp}.\qed

Table~\ref{tab:bipsimpevs} shows the number of bipartite graphs on $n$ vertices whose mixing matrices have given ranks.

\begin{table}[htbp]
    \centering
    \begin{tabular}{|c|c|c|c|}
    \hline \rule{0pt}{2.8ex} 
    $n$ & $\rk(\mxm)$ & \# bipartite graphs & \# simple eigenvalues\\
    \hhline{|=|=|=|=|}
    3&  2&  1&  1\\
    \hline
    3&  3&  0&  0\\
    \hhline{|=|=|=|=|}
    4&  2&  2&  1\\
    \hline
    4&  3&	0&  0\\
    \hline
    4&	4&	1&  0\\
    \hhline{|=|=|=|=|}
    5&	3&	3&  3\\
    \hline
    5&	4&	1&  0\\
    \hline
    5&	5&	1&  0\\
    \hhline{|=|=|=|=|}
    6&	2&	1&  1\\
    \hline
    6&	3&	6&  3\\
    \hline
    6&	4&	4&  0\\
    \hline
    6&	5&	4&  0\\
    \hline
    6&	6&	2&  0\\
    \hhline{|=|=|=|=|}
    7&	3&	0&  0\\
    \hline
    7&	4&	23& 20\\
    \hline
    7&	5&	3&  0\\
    \hline
    7&	6&	1&  0\\
    \hline
    7&	7&	3&  0\\
    \hhline{|=|=|=|=|}
    8&	2&	1&  1\\
    \hline
    8&	3&	5&  2\\
    \hline
    8&	4&	43& 24\\
    \hline
    8&	5&	51& 0\\
    \hline
    8&	6&	50& 0\\
    \hline
    8&	7&	21& 0\\
    \hline
    8&	8&	11& 0\\
    \hline
    \end{tabular}
    \caption{Number of bipartite graphs with given $\rk(\mxm)$.}
    \label{tab:bipsimpevs}
\end{table}

\begin{lemma}
	Let $S$ be a proper subset of the vertices of the graph $X$, and assume that for
	each vertex $a$ in $S$, there is an automorphism of $X$ with $a$ as its
	only fixed point. Then $\rk(\mxm)\ge|S|+1$.
\end{lemma}

\proof
We identify the automorphism group of $X$ with the set of permutation matrices
that lie in $\cmm(A)$. If $a\in S$, let $P_a$ be the automorphism of $X$ with
$a$ as its only fixed point. From Lemma~\ref{lem:decomp}, for each $P_a$, there is a matrix in $\Psi(\cD)$ with the same diagonal. These matrices and the identity matrix form a set of linearly independent matrices in $\Psi(\cD)$.
\qed

If $X$ is a Cayley graph for an abelian group
of odd order, then since the map that sends each group element to its inverse
gives rise to an automorphism of the Cayley graph with $1$ as its only fixed point,
it follows that each vertex is the unique fixed point of an automorphism of $X$.
This implies $\mxm$ must be invertible in this case. One easy consequence is that
the average mixing matrix of a cycle is invertible (although this is also an
a consequence of results in \cite{GodsilAverageMixing}).

Some of our numerical data indicates that for most graphs on $n$ vertices with
simple eigenvalues, the average mixing matrix has rank $n-1$. In view of Theorem 
\ref{lem:decomp}, this implies that, for these graphs, $\dim \cA_0 = 1$. Therefore 
any matrix that commutes with $A$ and has zero diagonal must be a scalar multiple of 
$A$. The following corollary is an immediate consequence.

\begin{corollary}
	Suppose $X$ is a connected graph with at least three vertices. If all 
	eigenvalues of $X$ are simple and $\rk(\mxm)=n-1$, then any automorphism 
	of $X$ has fixed points.\qed
\end{corollary}

\section{Open problems}

As the theory of the average mixing matrix is relatively new, there are many interesting problems one can ask. We discuss a few below.

One interesting question about the average mixing matrix concerns the non-negative rank of the average mixing matrix. The \textsl{non-negative rank} of a non-negative $n\times n$ matrix $A$ is the least number $k$, such that there are $k$  matrices $\{M_r\}_{r=1}^k$ of rank 1 with non-negative entries, such that $A = \sum_{r=1}^k M_r$.
If $X$ is a graph on $n >1$ vertices with simple eigenvalues, Lemma \ref{lem:simps} gives  that $\mxm$ has rank at most $n-1$. One can ask when the non-negative rank of $\mxm$ is equal to $n$. 

We may interpret the $(v,v)$ diagonal entry of $\mxm$ as the average probability of measuring at vertex $v$, after starting at vertex $v$. Thus questions about the trace of $\mxm$ are very natural. In particular, it is interesting to ask how graph invariants correspond to the trace of $\mxm$. Following in the vein of the questions about rank and trace, one can ask if the spectrum of $\mxm$ determines any graph properties. This gives rise to many natural questions about the average mixing matrix.

\end{document}